\documentclass[10pt,journal]{IEEEtran}
\usepackage{graphicx}          
\usepackage{amsmath}           
\usepackage{amssymb}
\usepackage{mathrsfs}
\usepackage{algorithm}
\usepackage{algpseudocode}
\usepackage{hyperref}
\usepackage[latin1]{inputenc}  

\usepackage{amsthm}

\theoremstyle{definition}

\theoremstyle{remark}

\theoremstyle{plain}

\DeclareMathOperator{\E}{E}
\DeclareMathOperator*{\argmin}{arg\,min}

\newcommand{\mbf}[1]{\mathbf{#1}}
\newcommand{\mbs}[1]{\boldsymbol{#1}}
\newcommand{\what}[1]{\widehat{#1}}
\newcommand{\wtilde}[1]{\widetilde{#1}}

\begin{document}

\title{Utilization of Noise-Only Samples in \\ Array Processing with Prior Knowledge}
\author{Dave Zachariah, Magnus Jansson and Mats Bengtsson\thanks{Copyright (c) 2013 IEEE. Personal use of this material is permitted. However, permission to use this material for any other purposes must be obtained from the IEEE by sending a request to pubs-permissions@ieee.org. The authors are with the ACCESS Linnaeus Centre, KTH Royal Institute of
Technology, Stockholm. E-mail:
$\{$dave.zachariah, magnus.jansson, mats.bengtsson$\}$@ee.kth.se. The research leading to these results has received funding from the European Research Council under the European Community's Seventh Framework Programme (FP7/2007-2013) / ERC grant agreement n$^\circ$ 228044 and the Swedish Research Council under contracts 621-2011-5847 and 621-2012-4134.}}

\maketitle

\begin{abstract}
For array processing, we consider the problem of estimating signals of interest, and their
directions of arrival (DOA), in unknown colored noise fields. We
develop an estimator that efficiently utilizes a set of noise-only samples and,
further, can incorporate prior knowledge of the DOAs with varying
degrees of certainty. The estimator is compared with state of the art
estimators that utilize noise-only samples, and the Cramér-Rao bound,
exhibiting improved performance for smaller sample sets and in poor signal conditions.
\end{abstract}


\section{Introduction}

Array signal processing has a wide range of
applications, including radar, communications, sonar, localization and
medical diagnosis \cite{VanTrees2002}. A central problem is that of direction of arrival (DOA) estimation. Several standard DOA estimators model the noise field as spatially white. When this assumption is violated by some arbitrary noise field, the performance can be severely degraded \cite{Viberg1993}. One option is to assume a more complex parametric noise model, cf. \cite{Friedlander&Weiss1995,Pesavento&Gershman2001}. Another option is to first use noise-only samples to estimate the noise statistics, then pre-whiten the subsequent data. This approach was, however, shown to be suboptimal \cite{Werner&Jansson2006}. Instead, \cite{Werner&Jansson2007} developed an approximate maximum likelihood (AML) estimator that uses the noise-only samples more efficiently. For further references to array processing in colored noise fields, cf. \cite{Werner&Jansson2006,Werner&Jansson2007}.

In certain scenarios, the DOA of the signals of interest are subject
to varying degrees of prior knowledge. State of the art methods
that incorporate such knowledge assume, however, that the noise field
is spatially white, cf. \cite{LinebargerEtAl1995,SteinwandtEtAl2011,BouleuxEtAl2009,Wirfalt&Jansson2012}.

In this paper we develop an estimator that is capable of utilizing the
noise-only samples more efficiently than the AML while also able to
incorporate prior knowledge of the DOAs of varying degrees of
certainty. The estimator is based on the maximum a posteriori (MAP)
framework and compared numerically with two state of the art
estimators and the Cramér-Rao bounds.

\emph{Notation:} $\mbf{A}^*$ and $\mathscr{C}(\mbf{A})$ denote the
Hermitian transpose and column space of $\mbf{A}$, respectively. The weighted inner product $\langle \mbf{x}, \mbf{y} \rangle_{\mbf{W}} \triangleq \mbf{y}^* \mbf{W} \mbf{x}$, where $\mbf{W} \succ \mbf{0}$ is positive definite.

\section{Problem formulation}

The output from an $m$-dimensional array is modeled by,
\begin{equation}
\mbf{y}(t) =
\begin{cases}
\mbf{n}(t), & t=-M+1, \dots, 0 \\
\mbf{A}(\mbs{\theta}) \mbf{s}(t) + \mbf{n}(t), & t=1, \dots, N,
\end{cases}
\end{equation}
where the noise plus interference is identically and independently complex Gaussian
distributed, $\mbf{n}(t) \sim \mathcal{CN}(\mbf{0}, \mbf{Q})$. Unlike, e.g., \cite{Friedlander&Weiss1995,Pesavento&Gershman2001} no structure of $\mbf{Q} \succ \mbf{0}$ is assumed. Rather the availability of a number of noise-only samples $M \geq m$ is assumed. During $N$
samples a $d$-dimensional signal $\mbf{s}(t) \in \mathbb{C}^d$ is
received, residing in a subspace parameterized by
$\mbf{A}(\mbs{\theta}) = [\mbf{a}(\theta_1) \cdots \mbf{a}(\theta_d) ]
\in \mathbb{C}^{m \times d}$. This models a set of $d$ narrowband
signals impinging on the array with directions of arrival (DOA)
$\mbs{\theta} = [\theta_1 \cdots \theta_d ]^\top$. The array is
assumed to be unambiguous, i.e., the columns of $\mbf{A}$ are linearly independent as long as $\theta_i \neq \theta_j$.

Let $\bar{\mbf{Y}} \triangleq [ \mbf{y}(-M+1)
\cdots \mbf{y}(0) ]$ and $\mbf{Y} \triangleq [ \mbf{y}(1) \cdots
\mbf{y}(N) ]$ denote the set of samples. The goal is to estimate $\mbs{\theta}$,
$\mbf{Q}$ and $\mbf{s}(t)$ from $\bar{\mbf{Y}}$ and $\mbf{Y}$. It is assumed that the problem is identifiable, however, the general conditions for this are difficult to derive.

The parameters are treated probabilistically. No prior knowledge is assumed about the signal $\mbf{s}(t)$ and covariance matrix $\mbf{Q}$, which are modeled by noninformative priors, $p(\mbf{Q}) \propto |\mbf{Q}|^{-(m+1)}$ and
$p(\mbf{S}) \propto 1$ \cite{Tiao&Zellner1964}, where $\mbf{S} = [ \mbf{s}(1) \cdots \mbf{s}(N) ]$, cf. \cite{Press2005} for further discussion on noninformative priors. Prior knowledge of
$\mbs{\theta}$ is modeled by independent von Mises distributions,
$\theta_i \sim \mathcal{M}(\mu_i, \kappa_i)$, which can be thought of
as a periodic analogue of the Gaussian distribution, where $\mu_i$ is the
expected value and $\kappa_i$ is a concentration parameter. When
$\kappa_i \rightarrow \infty$ it converges to a Gaussian
distribution with variance $1/\kappa_i$; when $\kappa_i =0$ it
corresponds to a noninformative prior, cf. \cite{EvansEtAl2000} and \cite{AbdiEtAl2002} for an illustration.

\vspace{-0.3cm}

\section{MAP estimator}

The MAP estimates of $\mbs{\theta}$, $\mbf{Q}$ and $\mbf{s}(t)$ are given by maximizing the posterior pdf $p( \mbs{\theta}, \mbf{Q}, \mbf{S} | \bar{\mbf{Y}}, \mbf{Y} )$. Equivalently, they are obtained by solving
\begin{equation}
  \max_{\mbs{\theta} \in  \Theta, \: \mbf{Q} \succ \mbf{0}, \: \mbf{S} \in \mathbb{C}^{d \times N} } J( \mbs{\theta}, \mbf{Q}, \mbf{S} ),
\label{eq:maxproblem}
\end{equation}
where $J( \mbs{\theta}, \mbf{Q}, \mbf{S} ) = J_1(\mbs{\theta}, \mbf{Q}, \mbf{S}) + J_2( \mbs{\theta} )$, and $J_1(\mbs{\theta}, \mbf{Q}, \mbf{S}) = \ln p(\bar{\mbf{Y}},  \mbf{Y} | \mbs{\theta},
\mbf{Q}, \mbf{S}  ) + \ln p(\mbf{Q}) + \ln p(\mbf{S})$ and $J_2(\mbs{\theta}) = \ln p( \mbs{\theta} )$ using Bayes' Rule.

\vspace{-0.3cm}

\subsection{Concentrated cost function}

First, we can simplify $J_1$ by noting the conditional independence
$p( \bar{\mbf{Y}}, \mbf{Y} | \mbs{\theta}, \mbf{Q}, \mbf{S}  ) = p( \bar{\mbf{Y}} | \mbf{Q}
)  p(\mbf{Y} | \mbs{\theta}, \mbf{Q}, \mbf{S}  ) $. For notational simplicity, let $\mbf{Q}_0 \triangleq \bar{\mbf{Y}}\bar{\mbf{Y}}^* / M$ and $\mbf{R}_0 \triangleq \mbf{Y}\mbf{Y}^* / N$ denote the sample covariance matrices. Define $\wtilde{\mbf{Y}} \triangleq \mbf{Y} - \mbf{A}\mbf{S}$,
and $\gamma \triangleq (M+N+m+1)$, so that
\begin{equation*}
\begin{split}
J_1 &= -M \ln | \mbf{Q} | - \text{tr}\{ \bar{\mbf{Y}} \bar{\mbf{Y}}^*
\mbf{Q}^{-1} \}\\
&\quad -N \ln | \mbf{Q} | - \text{tr}\{ \wtilde{\mbf{Y}} \wtilde{\mbf{Y}}^* \mbf{Q}^{-1} \}  - (m+1) \ln | \mbf{Q} | + K_1\\
&= -\gamma \ln |\mbf{Q}| - \text{tr} \{ ( M\mbf{Q}_0 + \wtilde{\mbf{Y}}
\wtilde{\mbf{Y}}^* ) \mbf{Q}^{-1}  \} + K_1,
\end{split}
\end{equation*}
where $K_1$ is a simple constant that can be omitted. The maximizing covariance matrix of $J_1(\mbs{\theta},\mbf{Q},\mbf{S})$ equals $\what{\mbf{Q}} = \frac{1}{\gamma} ( M \mbf{Q}_0 + \wtilde{\mbf{Y}} \wtilde{\mbf{Y}}^* )$ \cite{Stoica&Nehorai1995}. Then the concentrated cost function equals
\begin{equation*}
\begin{split}
J_1 &= -\gamma \ln |M\mbf{Q}_0 + (\mbf{Y} - \mbf{A}\mbf{S}) (\mbf{Y} -
\mbf{A}\mbf{S})^* |+ K'_1 \\
&= -\gamma \ln | \mbf{I}_m + M^{-1} \mbf{Q}^{-1}_0 (\mbf{Y} - \mbf{A}\mbf{S}) (\mbf{Y} -
\mbf{A}\mbf{S})^*  |+ K''_1 \\
&= -\gamma \ln |  \mbf{I}_N  +  M^{-1} (\mbf{Y} -
\mbf{A}\mbf{S})^* \mbf{Q}^{-1}_0 (\mbf{Y} - \mbf{A}\mbf{S})|+ K''_1,
\end{split}
\end{equation*}
where we used Sylvester's determinant theorem and $\mbf{Q}_0$ is invertible w.p.1.
Thus the inner argument of $J_1(\mbs{\theta}, \what{\mbf{Q}},
\mbf{S})$ is quadratic with respect to $\mbf{S}$. Since $-\ln | \cdot
|$ is a monotonically decreasing function on the set of positive
definite matrices, the stationary point is given at $\what{\mbf{S}} =
(\mbf{A}^*\mbf{Q}^{-1}_0 \mbf{A})^{-1} \mbf{A}^* \mbf{Q}^{-1}_0 \mbf{Y}$. As expected $\mbf{Q}^{-1}_0$ appears as a pre-whitening matrix.

Define $\mbs{\Phi}_{\mbf{A}} \triangleq \mbf{A}
(\mbf{A}^*\mbf{Q}^{-1}_0 \mbf{A})^{-1} \mbf{A}^* \mbf{Q}^{-1}_0$. This
matrix is the orthogonal projector onto $\mathscr{C}(\mbf{A})$ with
respect to the inner product $\langle \mbf{x} , \mbf{y}
\rangle_{\mbf{Q}^{-1}_0}$. Hence $\mbs{\Phi}^2_{\mbf{A}} =
\mbs{\Phi}_{\mbf{A}}$, $\mbf{Q}^{-1}_{0}\mbs{\Phi}_{\mbf{A}}$ is
Hermitian and $\mbs{\Phi}^\perp_{\mbf{A}} = \mbf{I}_m -
\mbs{\Phi}_{\mbf{A}}$ is the orthogonal projector onto $\mathscr{C}(\mbf{A})^\perp$ \cite{Rao&Rao1998}. Inserting the maximizer $\what{\mbf{S}}$ yields
\begin{equation}
\begin{split}
J_1 &= -\gamma \ln | \mbf{I}_N +  M^{-1} \mbf{Y}^* (\mbs{\Phi}^\perp_{\mbf{A}})^* \mbf{Q}^{-1}_0
\mbs{\Phi}^\perp_{\mbf{A}} \mbf{Y}  |+ K''_1 \\
&= -\gamma \ln | \mbf{I}_N + M^{-1} \mbf{Y}^* \mbf{Q}^{-1}_0 \mbs{\Phi}^\perp_{\mbf{A}} \mbs{\Phi}^\perp_{\mbf{A}} \mbf{Y} |+ K''_1 \\
&= -\gamma \ln | \mbf{I}_m + \alpha \mbf{Q}^{-1}_0 \mbs{\Phi}^\perp_{\mbf{A}} \mbf{R}_0 |+ K''_1, \\
\end{split}
\label{eq:J1}
\end{equation}
where $\alpha \triangleq N/M$. Next, the von Mises distribution yields \cite{EvansEtAl2000}
\begin{equation}
\begin{split}
J_2 &= \ln p(\mbs{\theta}) = \sum^d_{i=1} \kappa_i \cos(\theta_i - \mu_i) + K_2,
\end{split}
\label{eq:J2_prior}
\end{equation}
where $K_2$ is a constant. Finally, by combining \eqref{eq:J1} and \eqref{eq:J2_prior} the maximization problem \eqref{eq:maxproblem} can be recast as the concentrated minimization problem,
\begin{equation}
\hat{\mbs{\theta}} = \argmin_{\mbs{\theta} \in \Theta} \: \ln \left |\mbf{I}_m + \alpha \mbf{Q}^{-1}_0 \mbs{\Phi}^\perp_{\mbf{A}}(\mbs{\theta}) \mbf{R}_0 \right| + \varphi(\mbs{\theta}),
\label{eq:thetaMAP}
\end{equation}
where $\varphi (\mbs{\theta}) = - \sum^d_{i=1} \kappa_i \cos(\theta_i
- \mu_i) / \gamma$. This problem is nonconvex and a $d$-dimensional grid search may render it intractable.

\subsection{Iterative solution}

To make the problem tractable we exploit the decomposition property of
orthogonal projection matrices. For notational simplicity, let the $i$th column of $\mbf{A}$ be denoted as $\mbf{a}_i \in \mathbb{C}^{m \times 1}$ and the remaining columns $\mbf{A}_i \in \mathbb{C}^{m \times (d-1)}$. The projection operator can be decomposed as $\mbs{\Phi}_{\mbf{A}} = \mbs{\Phi}_{\mbf{A}_i} +
\mbs{\Phi}_{\tilde{\mbf{a}}_i}$, where $\tilde{\mbf{a}}_i =
\mbs{\Phi}^\perp_{\mbf{A}_i} \mbf{a}_i$. Further, the projector in \eqref{eq:thetaMAP} can be written as $\mbs{\Phi}^\perp_{\mbf{A}} = \mbs{\Phi}^\perp_{\mbf{A}_i} -
\mbs{\Phi}_{\tilde{\mbf{a}}_i}$, where
\begin{equation}
\begin{split}
\mbs{\Phi}_{\tilde{\mbf{a}}_i} &= \tilde{\mbf{a}}_i ( \tilde{\mbf{a}}^*_i \mbf{Q}^{-1}_0 \tilde{\mbf{a}}_i )^{-1} \tilde{\mbf{a}}^*_i \mbf{Q}^{-1}_0 \\
&= \mbs{\Phi}^\perp_{\mbf{A}_i} \mbf{a}_i \left( \mbf{a}^*_i   \mbf{Q}^{-1}_0 \mbs{\Phi}^\perp_{\mbf{A}_i} \mbs{\Phi}^\perp_{\mbf{A}_i} \mbf{a}_i \right)^{-1} \mbf{a}^*_i  \mbf{Q}^{-1}_0 \mbs{\Phi}^\perp_{\mbf{A}_i} \\
&= \frac{\mbs{\Phi}^\perp_{\mbf{A}_i} \mbf{a}_i \mbf{a}^*_i  \mbf{G}_i}{ \mbf{a}^*_i \mbf{G}_i \mbf{a}_i },
\end{split}
\label{eq:a_tilde}
\end{equation}
and where we defined $\mbf{G}_i \triangleq \mbf{Q}^{-1}_0 \mbs{\Phi}^\perp_{\mbf{A}_i}$ for notational simplicity. Then by defining $\mbs{\Psi}_i \triangleq \mbf{G}_i \mbf{R}_0 (\mbf{I}_m + \alpha \mbf{G}_i \mbf{R}_0 )^{-1} \mbf{G}_i$ and using \eqref{eq:a_tilde}, the determinant in \eqref{eq:thetaMAP} can be expressed as
\begin{equation*}
\begin{split}
 \left |\mbf{I}_m + \alpha \mbf{Q}^{-1}_0 \mbs{\Phi}^\perp_{\mbf{A}}\mbf{R}_0 \right| &=  \left | \mbf{I}_m + \alpha \mbf{G}_i \mbf{R}_0 - \alpha \mbf{Q}^{-1}_0 \mbs{\Phi}_{\tilde{\mbf{a}}_i} \mbf{R}_0 \right|\\
 &= \left | \mbf{I}_m + \alpha \mbf{G}_i \mbf{R}_0 \right | \left( 1 -  \alpha \frac{\mbf{a}^*_i \mbs{\Psi}_i \mbf{a}_i}{ \mbf{a}^*_i \mbf{G}_i \mbf{a}_i}  \right),
\end{split}
\end{equation*}
using the determinant theorem.

Following the alternating projections method in
\cite{Ziskind&Wax1988}, we can then relax \eqref{eq:thetaMAP} by
cyclicly optimizing over angle $\theta_i$ while keeping the
remaining angle estimates fixed in the vector $\mbs{\theta}'_i$ \cite{Stoica&Selen2004}.
This entails performing a series of one-dimensional grid searches
\begin{equation}
\hat{\theta}_i = \argmin_{\theta \in \Theta_i} \:  V( \theta ;
\mbs{\theta}'_i ),
\label{eq:V_search}
\end{equation}
where
\begin{equation}
V( \theta ; \mbs{\theta}'_i ) \triangleq \ln  \left( 1 -
  \alpha \frac{\mbf{a}^*(\theta) \mbs{\Psi}_i \mbf{a}(\theta)}{ \mbf{a}^*(\theta) \mbf{G}_i \mbf{a}(\theta)}  \right) + \varphi_i(\theta),
\label{eq:V_ap}
\end{equation}
and $\varphi_i(\theta) = -\kappa_i \cos(\theta - \mu_i)/\gamma$ for $i =
1,\dots, d$. The sequential search over a grid $\Theta_i$ of $g$ points is repeated until the difference between iterates, $|\Delta \hat{\theta}_i|$, is less than some tolerance.

For initialization we follow \cite{Ziskind&Wax1988}, starting with $\hat{\mbs{\theta}} = \varnothing$ and the angles $i=1,\dots,d$ sorted
with respect to $\kappa_i$ in descending order. This reduces the initial error in the search that arises when holding $\mbs{\theta}'_i$ constant. Initially,  $\Theta_i$ is  $[-90^\circ, 90^\circ]$ but the interval is subsequently refined in $L$ steps. The estimator is summarized in Algorithm~\ref{alg:MAP}. In the following, $\Theta_i$ is refined by reducing the interval by a half at each step and $\varepsilon_\ell$ is set to be equivalent of 2 grid points.

\begin{algorithm}
\caption{Alternating projections-based MAP estimator} \label{alg:MAP}
\begin{algorithmic}[1]
\State Input: $\bar{\mbf{Y}}, \mbf{Y}, \{ \mu_i, \kappa_i \}^d_{i=1}$ and $L$
\State Form $\mbf{Q}_0$, $\mbf{R}_0$, $\Theta^1_i$ and initialize $\hat{\mbs{\theta}} = \varnothing$
\For{$\ell = 1,\dots, L$}
    \Repeat
        \State For $i = 1, \dots, d$
        \State Form $\mbf{G}_i$ and $\mbs{\Psi}_i$
        \State $\hat{\theta}_i = \argmin_{\theta \in
          \Theta^\ell_i } V(\theta ; \mbs{\theta}'_i )$ using \eqref{eq:V_ap}
    \Until{ $|\Delta \hat{\theta}_i| < \varepsilon_\ell$ }
    \State Refine $\Theta^{\ell+1}_i, \forall i$
\EndFor
\State $\what{\mbf{S}} = \left(\mbf{A}^*(\hat{\mbs{\theta}})\mbf{Q}^{-1}_0 \mbf{A}(\hat{\mbs{\theta}}) \right)^{-1} \mbf{A}^*(\hat{\mbs{\theta}}) \mbf{Q}^{-1}_0 \mbf{Y}$
\State $\what{\mbf{Q}} = \left( M \mbf{Q}_0 + (\mbf{Y} -
\mbf{A}(\hat{\mbs{\theta}}) \what{\mbf{S}})( \mbf{Y} -
\mbf{A}(\hat{\mbs{\theta}}) \what{\mbf{S}}) \right)^* / \gamma$
\State Output: $\hat{\mbs{\theta}}$, $\what{\mbf{S}}$ and $\what{\mbf{Q}}$
\end{algorithmic}
\end{algorithm}

\section{Cramér-Rao bounds}

If the signals of interest are independent and identically distributed (i.i.d.) zero-mean Gaussian, i.e.,
$\mbf{s}(t) \sim \mathcal{CN}(\mbf{0}, \mbf{P})$, a Cramér-Rao bound
(CRB) for conditionally unbiased DOA estimators is given in \cite{Werner&Jansson2006}. The posterior CRB for random $\theta_i$ does not exist due to the circular von Mises distribution \cite{Routtenberg&Tabrikian2012}, but following \cite{WahlbergEtAl1991} we can formulate an approximate hybrid Cramér-Rao bound (ACRB) when the variance of random $\theta_i$ is small, using the result of \cite{Werner&Jansson2006}. The bound is given by
\begin{equation}
\mbf{C}_{\theta} =  \left( 2N\text{Re}\left\{ \mbf{D}^*( \mbs{\Gamma}^\top \otimes \mbf{Z}\mbs{\Pi}^\perp_{\mbf{ZA}}\mbf{Z} )\mbf{D} \right\} + \mbs{\Lambda}_\theta \right)^{-1},
\end{equation}
where $\mbf{Z}$ is the Hermitian square-root $\mbf{ZZ} = \mbf{Q}^{-1}$, $\mbf{D} = [\text{vec}(\partial_{\theta_1} \mbf{A}) \cdots \text{vec}(\partial_{\theta_d} \mbf{A}) ]$, $\mbs{\Pi}^\perp_{\mbf{ZA}}$ is the orthogonal projector onto $\mathscr{C}(\mbf{ZA})^\perp$ and $\mbs{\Gamma} = \mbf{P}\mbf{A}^* \mbf{Z}^* \mbf{E}_s(\mbs{\Lambda}_s + \alpha \mbf{I}_d)^{-1}\mbf{E}^*_s \mbf{Z}\mbf{A} \mbf{P}$. Here $\mbf{E}_s$  and $\mbs{\Lambda}_s$ are given by the eigendecomposition of $\mbf{ZRZ}$. See \cite{Werner&Jansson2006} for details. The matrices dependent on $\theta_i$ are evaluated at the expected values $\mbf{\mu}_i$. Finally, $\mbs{\Lambda}_\theta = \text{diag}(\lambda_1, \dots, \lambda_d)$ embodies the prior information. The diagonal elements $\lambda_i$ equal $\kappa_i$ or $0$ depending on whether $\theta_i$ is treated as a random or deterministic quantity, respectively.

\section{Experimental results}

We consider a uniform linear array (ULA) with half-wave length
separation. For comparison, we also consider two state of the art
estimators: The optimally weighted MODE estimator, denoted W-MODE \cite{Werner&Jansson2006}, and the approximate maximum likelihood estimator, denoted C-MODE \cite{Werner&Jansson2007}. Both are asymptotically efficient. We evaluate the estimators using the root mean square error $\text{RMSE}(\hat{\theta}_i) \triangleq \sqrt{ \E[ \tilde{\theta}^2_i ] }$, where $\tilde{\theta}_i$ is the estimation error. The RMSE is evaluated numerically using $5\cdot 10^3$ Monte Carlo runs.

\subsection{Setup}

We consider $d=3$ correlated Gaussian signals $\mbf{s}(t)$ with covariance
matrix $\mbf{P} = \mbf{I}_d + \rho \mbf{T} + \rho^* \mbf{T}^*$, where $0 \leq
|\rho| < 1$ and $\mbf{T}$ is a strictly lower-triangular matrix with nonzero elements
equal to 1. The first angle, $\theta_1$, is considered with an expected value $\mu_1$ and certainty given by $\kappa_1 = 10^5$, corresponding to a standard deviation of about
$0.18^\circ$, while the prior knowledge of the remaining angles, $\theta_2$ and $\theta_3$, is noninformative, i.e.,
$\kappa_2 = \kappa_3 = 0$. Then the first DOA, $\theta_1$, is
randomized according to $\mathcal{M}(\mu_1, \kappa_1)$ with $\mu_1 = -35^\circ$ \cite{Berens2009}, while the remaining DOAs are fixed as $\theta_2 = 15^\circ$ and $\theta_3 = 20^\circ$.

The unknown noise field is modeled as spatially correlated noise plus $\tilde{d} =
3$ interferers with DOAs $\tilde{\mbs{\theta}} = [-40^\circ, \: -10^\circ, \:
40^\circ]^\top$. The noise covariance matrix is $\mbf{Q} = \mbf{Q}' +
\mbf{A}(\tilde{\mbs{\theta}})\wtilde{\mbf{P}}
\mbf{A}^*(\tilde{\mbs{\theta}})$. Here $\{ \mbf{Q}' \}_{ij} = \sigma^2
a^{|i-j|}$, where $a \in [0,1)$ controls the spatial correlation,
and $\wtilde{\mbf{P}} = \tilde{\sigma}^2 \mbf{I}_{\tilde{d}}$.

We consider an array of $m=10$ elements, with sample ratio $\alpha=1$ and
spatial signal and noise correlations $\rho = 0.9$ and $a = 0.5$,
respectively. Three parameters are varied: (a) the number of samples
$M$, (b) signal to noise ratio $\text{SNR} \triangleq \text{tr}\{
\mbf{P} \} / \text{tr}\{ \mbf{Q}' \} $ and (c) interference to noise
ratio $\text{INR} \triangleq \text{tr}\{ \wtilde{\mbf{P}} \}/ \text{tr}\{ \mbf{Q}' \} $.

For MAP, $\Theta_i$ is a grid of $g = 500$ points and the grid refinement is repeated $L=10$ times, yielding a resolution limit of $180 /(2^{L-1} g) \approx 7 \times 10^{-4}$ degrees. For C-MODE and W-MODE, we use 3 iterations as in \cite{Werner&Jansson2007}.

\subsection{Results}

\begin{figure}
  \begin{center}
    \includegraphics[width=0.8\columnwidth]{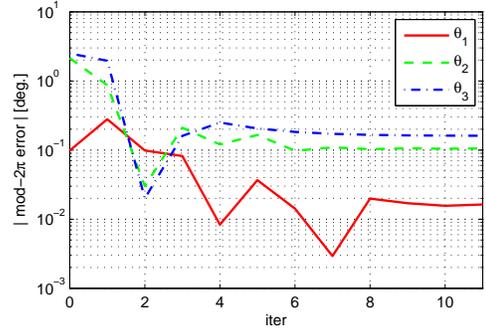}
  \end{center}
  \vspace{-0.35cm}
  \caption{Example convergence of the MAP estimates. Absolute error $|\tilde{\theta}_i|$ versus iteration for a typical realization with $M=100$, $\text{SNR} = 5$~dB and $\text{INR} = 5$~dB. Each iteration corresponds to $d$ grid searches \eqref{eq:V_search}. The algorithm terminated at the 11th iteration.}
  \label{fig:convergence}
\end{figure}

For the current nonoptimized implementations of the estimators, the execution time for a typical realization is 2, 61 and 550 milliseconds for W-MODE, C-MODE and MAP, respectively. This should be compared with 2 milliseconds required to compute the sample covariance matrices for $M=N=10^4$. Fig.~\ref{fig:convergence} illustrates the convergence behavior of MAP for a typical realization. Note that $\hat{\theta}_1$ starts with a lower error due to prior knowledge. At each iteration the cost in \eqref{eq:thetaMAP} declines.

\begin{figure}
  \begin{center}
    \includegraphics[width=0.88\columnwidth]{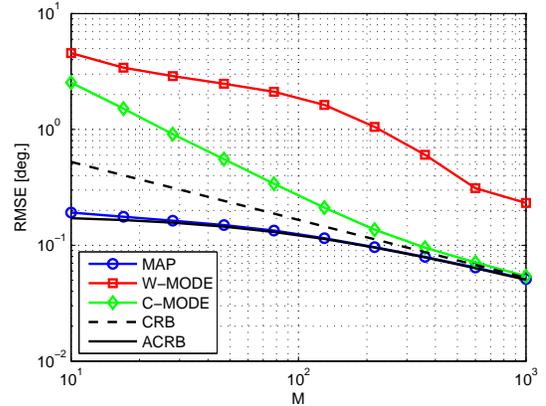}
  \end{center}
  \vspace{-0.3cm}
  \caption{$\text{RMSE}(\hat{\theta}_1)$ versus sample size $M$ at SNR=5~dB and INR=5~dB.}
  \label{fig:M_01}
\end{figure}

\begin{figure}
  \begin{center}
    \includegraphics[width=0.88\columnwidth]{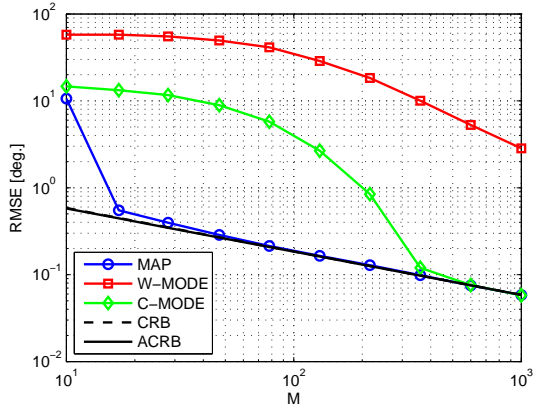}
  \end{center}
  \vspace{-0.3cm}
  \caption{$\text{RMSE}(\hat{\theta}_3)$ versus sample size $M$ at SNR=5~dB and INR=5~dB.}
  \label{fig:M_03}
\end{figure}

Figs.~\ref{fig:M_01} and \ref{fig:M_03} show the RMSE performance with increasing sample size $M$, for $\theta_1$ and $\theta_3$, respectively. The first angle, $\theta_1$, with an informative prior is  surrounded by interferers at $-40^\circ$ and $-10^\circ$. In this case the ACRB is visibly below CRB for low $M$. We see that MAP is able to improve on the prior knowledge of $\theta_1$. The DOAs with noninformative prior knowledge, $\theta_2$ and $\theta_3$, are closely spaced and surrounded by interferers at $-10^\circ$ and $40^\circ$. In this case CRB and ACRB are virtually identical. For both DOAs, MAP approaches the ACRB at low $M$ while the alternative estimators require more than an order of magnitude more samples to close the gap. Thus while MAP is more computationally complex than W-MODE and C-MODE, for a given performance level it can substantially reduce the number of snapshots to acquire and compute $\mbf{Q}_0$ and $\mbf{R}_0$. Further this enables less restrictive assumptions on the stationarity period of the noise.

\begin{figure}
  \begin{center}
    \includegraphics[width=0.88\columnwidth]{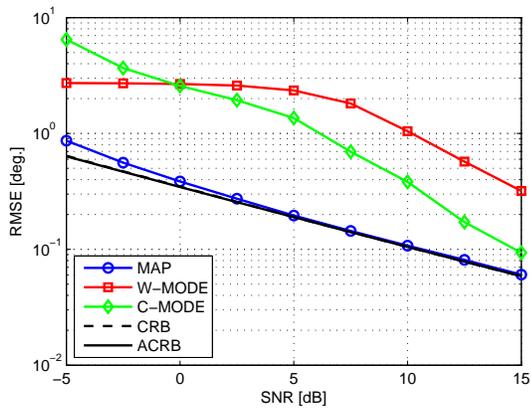}
  \end{center}
  \vspace{-0.3cm}
  \caption{$\text{RMSE}(\hat{\theta}_2)$ versus SNR at $M = 100$ and INR=5~dB.}
  \label{fig:SNR_02}
\end{figure}

\begin{figure}
  \begin{center}
    \includegraphics[width=0.88\columnwidth]{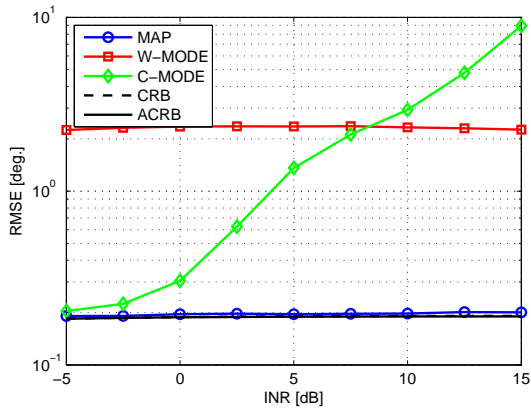}
  \end{center}
  \vspace{-0.3cm}
  \caption{$\text{RMSE}(\hat{\theta}_2)$ versus INR at $M = 100$ and SNR=5~dB.}
  \label{fig:INR_02}
\end{figure}

The RMSE performance for $\theta_2$ is shown in Fig.~\ref{fig:SNR_02} as a function of SNR. While the other estimators approach the CRB when the signal to noise ratio reaches 20~dB, i.e., substantially above INR=5~dB, MAP matches the average performance at lower SNR. The key explanation for the advantage of MAP over C-MODE and W-MODE is its resilience to interfering sources as illustrated in Fig.~\ref{fig:INR_02}. Unlike the other estimators, MAP forms an optimal estimate of the noise covariance matrix without approximations. This allows it to reject the noise even when INR is substantially greater than SNR.



\section{Conclusion}

We have developed a DOA and signal of interest estimator using the MAP framework, that utilizes noise only-samples and is capable of incorporating prior knowledge of the DOAs. By forming an optimal estimate of the noise covariance matrix, the DOA estimates are especially resilient to strong interferers. An alternating projections-based method was used to solve the resulting optimization problem. Finally, the resulting estimator was compared with the state of the art C-MODE and W-MODE as well as the Cramér-Rao bounds, exhibiting significantly improved average performance at smaller sample sets and deteriorating signal conditions.

\vspace{-0.27cm}

\bibliographystyle{ieeetr}
\bibliography{refs_array}

\end{document}